\documentclass[11pt]{amsart}
\usepackage{amsmath}
\usepackage{xypic}
\usepackage{amssymb}
\usepackage{hyperref}
\usepackage{amsthm}

\newtheorem{theorem}{Theorem}[section]
\newtheorem{proposition}[theorem]{Proposition}
\newtheorem{lemma}[theorem]{Lemma}
\newtheorem{definition}[theorem]{Definition}

\newtheorem{corollary}[theorem]{Corollary}
\newtheorem{remark}[theorem]{Remark}
\newtheorem{example}[theorem]{Example}

\begin{document}

\def\be{\begin{equation}}
\def\ee{\end{equation}}

\def\ra#1{\mathop{\longrightarrow}\limits^{#1}}
\def\rab#1{\mathop{\longrightarrow}\limits_{#1}}
\def\lab#1{\mathop{\longleftarrow}\limits_{#1}}
\def\la#1{\mathop{\longleftarrow}\limits^{#1}}
\def\da#1{\downarrow\rlap{$\vcenter{\hbox{$\scriptstyle#1$}}$}}
\def\ua#1{\uparrow\rlap{$\vcenter{\hbox{$\scriptstyle#1$}}$}}
\def\sea#1{\mathop{\searrow}\limits^{#1}}
\def\nea#1{\mathop{\nearrow}\limits^{#1}}
\def\swa#1{\mathop{\swarrow}\limits^{#1}}
\def\dirlim#1{\lim_{\textstyle{\mathop{\rightarrow}\limits_{#1}}}}
\def\invlim#1{\lim_{\textstyle{\mathop{\leftarrow}\limits_{#1}}}}

\def\AA{\mathbb{A}}
\def\DD{\mathbb{D}}
\def\QQ{\mathbb{Q}}
\def\ZZ{\mathbb{Z}}
\def\NN{\mathbb{N}}
\def\GG{\mathbb{G}}
\def\WW{\mathbb{W}}
\def\RR{\mathbb{R}}
\def\CC{\mathbb{C}}
\def\FF{\mathbb{F}}
\def\HH{\mathbb{H}}
\def\LL{\mathbb{L}}

\def\qed{$\,\square$}

\def\tensor{\otimes}
\def\dsum{\oplus}
\def\Tensor{\mathop\bigotimes}
\def\Dsum{\mathop\bigoplus}
\def\img{\mathrm{im}}

\def\bR{\mathbf{R}}

\def\cW{\mathcal{W}}
\def\cC{\mathcal{C}}
\def\cF{\mathcal{F}}
\def\cG{\mathcal{G}}
\def\cO{\mathcal{O}}
\def\cL{\mathcal{L}}
\def\cM{\mathcal{M}}
\def\cN{\mathcal{N}}

\def\bSet{\mathbf{Set}}
\def\bMod{\mathbf{Mod}}
\def\bAb{\mathbf{Ab}}
\def\bC{\mathbf{C}}
\def\bA{\mathbf{A}}
\def\bMag{\mathbf{Mag}}
\def\bV{\mathbf{V}}
\def\bAss{\mathbf{Ass}}
\def\bCom{\mathbf{Com}}
\def\bAlt{\mathbf{Alt}}
\def\bAAlt{\mathbf{AAlt}}
\def\bDelta{\mathbf{\Delta}}
\def\bSigma{\mathbf{\Sigma}}
\def\bOrd{\mathbf{Ord}}
\def\bLie{\mathbf{Lie}}
\def\bNov{\mathbf{Nov}}
\def\bVin{\mathbf{Vin}}
\def\bAlg{\mathbf{Alg}}
\def\bCA{\mathbf{CAlg}}
\def\bLeib{\mathbf{Leib}}
\def\bLev{\mathbf{Lev}}
\def\bRMod{\mathbf{RMod}}
\def\bLMod{\mathbf{LMod}}

\def\Hom{\mathrm{Hom}}
\def\Map{\mathrm{Map}}
\def\c{\mathrm{c}}
\def\Ab{\mathrm{Ab}}
\def\End{\mathrm{End}}
\def\pr{\mathrm{pr}}
\def\Tens{\mathrm{Tens}}
\def\gr{\mathrm{gr}}
\def\Ext{\mathrm{Ext}}
\def\sgn{\mathrm{sgn}}
\def\res{\mathrm{res}}
\def\rank{\mathrm{rank}}
\def\Der{\mathrm{Der}}
\def\inj{\mathrm{in}}

\def\intersect{\mathop\bigcap}
\def\iso{\cong}
\def\cotensor{\square}

\def\bDPA{\mathbf{DPAlg}}

\def\ext{{e}}

\author{Aseel Kmail}
\address{Department of Mathematics, Arab American University
Jenin, Palestine}
\email{a.kmail45@student.aaup.edu}

\author{Julia Kozak}
\address{Department of Mathematics, Massachusetts Institute of Technology,
Cambridge, MA 02139}
\email{juliakzk@mit.edu}

\author{Haynes Miller}
\address{Department of Mathematics, Massachusetts Institute of Technology,
Cambridge, MA 02139}
\email{hrm@math.mit.edu}

\title
{Divided powers and K\"ahler differentials}


\begin{abstract}
Divided power algebras form an important variety of non-binary universal 
algebras. We identify the universal enveloping algebra and K\"ahler 
differentials associated to a divided power algebra over a general commutative
ring, simplifying and generalizing work of Roby and Dokas. 
\end{abstract}

\setcounter{equation}{0}

\maketitle

\section{Introduction}

Divided power algebras were introduced by Henri Cartan \cite{cartan} to 
describe the homology of Eilenberg Mac Lane spaces. They have subsequently
been intensively studied \cite{roby-pd} and have
played important roles in other parts of mathematics, such as algebraic
geometry, where they form the basis of the construction of crystalline
cohomology \cite{berthelot-ogus}. They constitute an important example of a
variety of algebras that are ``non-binary,'' in the sense that their 
structure is not entirely encoded in a binary product. 

In \cite{quillen}, Dan Quillen described a uniform process for defining a
cohomology theory in any one of a wide class of algebraic structures. 
An important role in that construction is played by 
the category of ``Beck modules'' (see \cite{barr}, for example)
for an algebra $A$ in a specified variety $\bV$. This is the category,
recognized long ago by Sammy Eilenberg \cite{eilenberg}
as providing the appropriate
meaning of a ``representation'' in a general context, consists of the 
abelian objects in the slice category $\bV/A$. In linear cases it can
be identified with the category of modules over a unital associative 
algebra $U(A)$, 
the ``universal enveloping algebra'' of $A$. Beck modules form the
coefficients in the Quillen cohomology theory defined on $\bV$, and in fact
Quillen homology is an appropriately defined derived functor of the
abelianization functor, evaluated on the terminal object $1_A:A\downarrow A$
in $\bV/A$. If $\bV$ is the variety of commutative rings, $\Ab_A(1_A)$
is the $A$-module of K\"ahler differentials, and this suggests defining
$\Omega^\bV_A=\Ab_A(1_A)$ in a general variety of algebras $\bV$.

This construction has been considered in detail by Ionnnis Dokas 
\cite{dokas-differentials,dokas-triple} in case
one is working with divided power algebras over a field. 
The goal of the present paper is to show
how this works out over a general ring. Dokas proves the interesting 
result that if $A$ is a $DP$ algebra over a field then the module of 
divided power K\"ahler differentials is simply a $DP$ Beck module 
structure on the usual commutative algebra module of K\"ahler differentials. 
We show that the same
result holds in general, and along the way we simplify some of his 
arguments. 

The classical module of K\"ahler differentials of divided power algebras was
also the subject of a study by Norbert Roby \cite{roby-differentials}. 
He described the K\"ahler differentials of a free divided power algebra.
We show how his result follows easily from the identification 
$\Omega^{DP}_A=\Omega^{CA}_A$. 

In the first section of this short paper we gather some reminders about
divided power algebras. We find it convenient to work with {\em non-unital}
commutative algebras. We then study the structure of abelian divided
power algebras -- those whose product is trivial. These coincide with 
Beck modules over the zero divided power algebra. In \S4 we identify
the universal enveloping algebra of a general divided power algebra.
Next we define the module of K\"ahler differentials, and state and prove 
the main theorem. Finally, in \S6, we give some consequences of this theorem.

\smallskip\noindent
{\bf Acknowledgements.}
This work was carried out under the auspices of the ``PalUROP'' program,
designed to foster
collaborative research projects involving students from MIT and Palestinian
universities.
The second author acknowledges support by the MIT UROP office.

\section{Divided power algebras and their universal enveloping algebras}

We recall the definition of divided power algebras and some of their basic
properties. We will always work over a fixed commutative ring $R$, 
and an un-decorated ``$\tensor$'' will indicate the tensor product over $R$. 
By an ``algebra'' we mean a commutative non-unital $R$-algebra. 

\begin{definition} A {\em $DP$ structure} on the algebra $A$ 
is a family of maps $\gamma_i:A\to A,i>0$ such that for all $a,b\in A$ and 
$r\in R$,
\begin{gather*}
\gamma_1(a)=a\\
\gamma_n(a+b)=\gamma_n(a)+\sum_{i+j=n}\gamma_i(a)\gamma_j(b)+\gamma_n(b)\\
\gamma_n(ab)=a^n\gamma_n(b)\quad\text{\rm and}\quad \gamma_n(rb)=r^n\gamma_n(b)\\
    \gamma_m(a)\gamma_n(a)=\frac{(m+n)!}{m!n!}\gamma_{m+n}(a)\\
    \gamma_m(\gamma_n(a))=\frac{(mn)!}{m!(n!)^m}\gamma_{mn}(a)\,.
\end{gather*}
A morphism of $DP$ algebras is an algebra map commuting with these
``divided power'' operations.
\end{definition}

Write $\bDPA_R$ for the category of $DP$ algebras.

\begin{remark} {\em The commutative ring $R$ might be ``large,'' including
the nonunital algebra $A$ as an ideal. 
Then one may define $\gamma_0=1$ and the axioms become 
slightly more compact; this is the expression one finds in 
\cite{berthelot-ogus} and elsewhere.

But $R$ may also be ``small.'' 
In fact any nonunital $R$-algebra $A$ is an ideal in
the augmented unital $R$-algebra 
\[
A_+=A\oplus R
\]
with product $(a,r)(b,s)=(ab,rb+sa)$.
This is the context discussed by Roby \cite{roby-pd} and Dokas 
\cite{dokas-triple,dokas-differentials}.
(Our use of the subscripted $+$ is the opposite of theirs, but agrees 
with the usage in \cite{alt-alg}.)
}\end{remark}

The following well-known facts will be useful. 

\begin{proposition}[\cite{roby-pd}]
The forgetful functor $\bDPA_R\to\bMod_R$ 
has a left adjoint, sending the $R$-module $V$ to the ``free $DP$ algebra'' 
$\Gamma_R(V)$ generated by $V$. 
\end{proposition}

The coproduct of two algebras, $A$ and $B$, is the algebra
\[
A\textstyle{\coprod} B=(A\tensor R)\oplus(A\tensor B)\oplus(R\tensor B)
\]
with the evident product. The canonical inclusions take $a\in A$ to $a\tensor1$ 
and $b\in B$ to $1\tensor b$. 

\begin{proposition}[{\cite{roby-pd}}]
Let $A$ and $B$ be $DP$ algebras. Then there is a unique
$DP$ structure on $A\coprod B$ such that the two inclusions are $DP$ maps,
and it serves as the coproduct in $\bDPA_R$.
\end{proposition}

The product of two $DP$ algebras is a 
$DP$ structure on their product as algebras, which has as its 
underlying $R$-module the product $R$-module. The terminal $DP$ algebra
is the unique structure on the trivial $R$-module $0$. Any algebra
has a unique ``point,'' i.e. a unique map from the terminal algebra. 

\section{Abelian $DP$-algebras}

We are interested in abelian group objects in the category $\bDPA_R$. 
A unital product on an $R$-algebra $A$ is an algebra map 
$\mu:A\times A\to A$ such that $\mu(a,0)=a$ and $\mu(0,b)=b$. 
From this we find that 
\[
ab=\mu(a,0)\mu(0,b)=\mu(0,0)=0\,.
\]
Also, $\mu$ is an {\em additive} map, so 
\[
\mu(a,b)=\mu((a,0)+(0,b))=\mu(a,0)+\mu(0,b)=a+b
\]
so the product is none other than the sum in the $R$-module $A$,
which is $R$-linear by distributivity. This product is thus automatically 
an abelian group structure, 
so the abelian objects in $\bCA_R$ are the algebras with trivial product.
Since the categorical product in $\bDPA_R$ is the same, an abelian $DP$ 
algebra also has 
trivial algebra structure. This imposes strong restrictions on the divided 
powers.

\begin{proposition} A $DP$ algebra is abelian as a $DP$ algebra if and
only if it is abelian as a commutative algebra. In such a $DP$ algebra,
$\gamma_n=0$ unless $n$
is a power of a prime. For any prime number $p$, $\gamma_p$ is additive,
$p\gamma_p=0$, $\gamma_{p^e}=\gamma_p^e$, and, for all $a\in A$ and $r\in R$,
$\gamma_p(ra)=r^p\gamma_p(a)$.
\end{proposition}
\begin{proof}
We begin by verifying these properties.
Suppose $A$ is a $DP$ algebra that is abelian as a commutative algebra.
The deviation from additivity of $\gamma_n$ is a sum of products, so 
vanishes if products vanish. 
If $i+j=n$ with $i,j\geq1$, the product rule for divided powers shows that 
$(i,j)\gamma_n=0$ (where $(i,j)=(i+j)!/i!j!$). Since \cite{ram}
\[
\gcd\{(i,j):i,j\geq1,i+j=n\}=
\begin{cases}
1 & \text{ if}\,\,n\,\,\text{is not a prime power}\\
p & \text{ if}\,\, n=p^e
\end{cases}
\]
we see that $\gamma_n=0$ if $n$ is not a prime power and $p\gamma_{p^e}=0$. 
Finally, following \cite[7-09]{cartan}, we claim that
\[
\gamma_{kp}\equiv \gamma_k\gamma_p \mod p\,.
\]
That is to say,
\[
\frac{(kp)!}{k!(p!)^k}\equiv1\mod p\,.
\]
To see this, write the numerator as
\[
\prod_{i=1}^p\big(i(p+i)(2p+i)\cdots((k-1)p+i)\big)\,.
\]
The denominator is $((p-1)!)^k$ times the $i=p$ factor in the numerator.
For $i<p$, the $i$th factor in the numerator
is congruent mod $p$ to $i^k$, so the product of those terms is congruent
mod $p$ to $((p-1)!)^k$. 

It follows by induction that 
$\gamma_{p^e}\equiv\gamma_p^e\mod p$.

Now observe that if $A$ is $DP$ algebra with trivial product then the
addition map, which is the only candidate for an abelian group product,
is indeed a $DP$ algebra map. This follows from linearity of the
divided powers on such an algebra.  
\end{proof}

Define the unital associative algebra
\[
U(0)=R[\phi_p:p \,\,\text{prime}]/(p\phi_p)\,,
\]
with the twisted product defined by $\phi_pr=r^p\phi_p$. 
Note that the map $R\to U(0)$ is not generally central; $U(0)$ is not an 
$R$-algebra, but it is an ``$R$-bimodule-algebra'' in the sense that it is
a monoid in the monoidal category of $R$-bimodules. We have proved that 
an abelian $DP$ algebra structure on $A$ is recorded by a left $U(0)$ action
$U(0)\tensor_RA\to A$. The $\phi$'s record the divided powers; we have
changed notation from $\gamma$ to $\phi$ for clarity of later constructions. 
If $p$ and $q$ 
are distinct primes, then both $p$ and $q$ kill $\phi_p\phi_q$, which 
therefore vanishes. As a notational convenience, let $\phi_{p^e}=\phi_p^e$ 
and $\phi_n=0$ if $n$ is not a prime power.  

Conversely, any such action determines an abelian
$DP$ algebra structure on $A$. So:

\begin{lemma}
The category of abelian $DP$ algebras is equivalent to the category of
left $U(0)$-modules
\end{lemma}

A $DP$ ideal in a $DP$ algebra $A$ is an ideal in $A$ that is closed under 
formation of divided powers. It is easy to see that if $I$ is a $DP$ ideal 
in $A$ 
then the divided powers descend to give a divided power structure on $A/I$,
their non-additivity notwithstanding. 

The ideal $A^2$ is closed under the divided powers, so $A/A^2$ is an example
of an abelian divided power algebra; indeed, it is the abelianization of $A$. 

\section{Beck modules}

Fix a $DP$ algebra $A$ and consider the slice category $\bDPA_R/A$. 
An object of $\bDPA_R/A$ is a $DP$ algebra $B$ equipped with a $DP$ map 
$\pi:B\downarrow A$. A morphism from $\pi':B'\downarrow A$ to 
$\pi:B\downarrow A$ is a $DP$ algebra map $f:B'\to B$ such that $\pi f=\pi'$.

\begin{definition} A {\em $DP$ $A$-module} is an abelian object in 
$\bDPA_R/A$. 
\end{definition}

These are the ``Beck modules'' in the theory of $DP$ algebras.
Write $\bMod_A$ for the category of $DP$ $A$-modules. The following
was observed by Dokas \cite[\S3]{dokas-triple} if $R$ is a field.

\begin{proposition} $\bMod_A$ is an abelian category, equivalent to the 
category of left modules over the $R$-bimodule-algebra $U(A)$ 
given by 
\[
U(A)=A_+\tensor_RU(0)
\]
as an $R$-bimodule, with product determined by 
\begin{gather*}(a\tensor1)(b\tensor1)=
ab\tensor 1\,,\quad (1\tensor u)(1\tensor v)=1\tensor uv\\
(a\tensor1)(1\tensor u)=a\tensor u\\
(1\tensor\phi_p)(a\tensor 1)=0
\end{gather*}
with $a,b\in A$ and $u,v\in U(0)$ and $p$ prime.
\end{proposition}
\begin{proof}
Given an abelian object $\pi:E\downarrow A$ over $A$ in $\bDPA_R$, let 
$M=\ker\pi$. This kernel is a sub $DP$ algebra of $E$. The abelian structure
of $E\downarrow A$ restricts to an abelian structure on $M\downarrow 0$,
so the product on $M$ is trivial and the divided powers are given by an
action of $U(0)$ on $M$. 

The unital algebra $A_+$ acts on $E$, and $M$ is a submodule. It remains
to describe how the action of $A_+$ interacts with the divided power 
structure. But for $n>1$, in $A\oplus M$,
\[
(0,\gamma_n(ax))=
\gamma_n(0,ax)=\gamma_n((a,0)(0,x))=\gamma_n(a,0)(0,x)^n=0
\]
so $\gamma_n(ax)=0$. 

Conversely, given a $U(A)$-module $M$, the $R$-module $A\oplus M$ becomes
an abelian object over $A$ with 
\[
(a,x)(b,y)=(ab,ay+bx)\,,\quad 
\gamma_n(a,x)=\Big(\gamma_na,\phi_nx+\sum_{i+j=n}\gamma_i(a)\phi_j(x)\Big)\,.
\]
These constructions are inverse to each other.
\end{proof}

It's interesting to observe that the $DP$ structure of $A$ plays no role in 
this description of the category of $DP$ $A$-modules. 

For example, fix a prime $p$ and suppose that $R$ is a $\ZZ_{(p)}$-algebra.
Then $U(A)$ is too, so $\phi_q=0$ for primes $q$ other
than $p$, and
\[
U(A)=A_+\tensor_R R[\phi_p]/(p\phi_p)
\]
with product twisted by $\phi_pr=r^p\phi_p$ for $r\in R$ and 
$\phi_pa=0$ for $a\in A$. Compare with \cite{dokas-differentials}.

\begin{remark}{\em Divided power algebras often occur in a graded setting;
this is how Cartan encountered them in \cite{cartan}, for example. 
If the Koszul sign rule applies, one may divide by the ideal generated by
the elements of odd degree to obtain a graded algebra that is commutative
in the sense that $ab=ba$. 

In this graded context, it is appropriate to grade $U(0)$ on the commutative
monoid $\ZZ_{>0}^\times$ of positive integers under multiplication;
so for example $|1|=|\phi_1|=1$. This monoid
acts additively on $\ZZ$, and the action of the divided powers is compatible
with that action. The algebra $U(A)$ is naturally $\ZZ$-graded; an element 
$a\tensor \phi_n\in A_+\tensor_R U(0)$ has degree $n|a|$. Then the category of
abelian objects in the category of graded $DP$ algebras over $A$ is equivalent
to the category of graded left $U(A)$-modules. 
}\end{remark}

\section{Derivations and K\"ahler differentials}

The category $\bDPA_R$ has its own proper theory of derivations and 
differentials. 

\begin{definition} Let $A$ be a $DP$ algebra and $M$ a $DP$ $A$-module. 
A $DP$ {\em derivation} is an $R$-linear map $s:A\to M$ such that 
$a\mapsto (a,s(a))$ is a $DP$ algebra section of $\pr_1:A\oplus M\downarrow A$.
\end{definition}

\begin{lemma} An $R$-linear map $s:A\to M$ is a $DP$ derivation if and only
if 
\begin{gather*}
s(ab)=as(b)+bs(a)\\
s(\gamma_n(a))=\phi_n(sa)+\sum_{i+j=n}\gamma_i(a)\phi_j(sa)\,.
\end{gather*}
\end{lemma}
\begin{proof}
The first equality follows from the map $a\mapsto(a,sa)$ being an algebra map: 
\[
(ab,s(ab))=(a,s(a))(b,s(b))=(ab,as(b)+bs(a))\,.
\]
The second follows from commutation with divided powers:
\begin{gather*}
(\gamma_na,s(\gamma_na))=\gamma_n(a,sa)=\gamma_n((a,0)+(0,sa))\\
=\gamma_n(a,0)+\sum_{i+j=n}(\gamma_ia,0)(0,\phi_j(sa))+\gamma_n(0,sa)\\
=(\gamma_na,0)+\sum_{i+j=n}(0,(\gamma_ia)(\phi_j(sa)))+(0,\phi_n(sa)) \,.
\end{gather*}
Equating second entries gives the result. 

Conversely, these conditions imply that $a\mapsto(a,s(a))$ is a $DP$ algebra
map splitting the projection. 
\end{proof}

For example, if $A$ is a $\ZZ_{(p)}$-algebra, a $DP$ derivation is an 
algebra derivation $s$ such that $s(\gamma_pa)=\phi_p(sa)$.

\begin{remark}\label{rem-implications}
{\em
We note some implications of these equations. First of all, 
\[
\phi_p(s(ab))=0
\]
since $s(ab)=as(b)+bs(a)$ and $\phi_p$ kills elements of the form $ax$. 
Next, 
\[
\phi_p(s(\gamma_qa))=\phi_p\phi_q(sa)+\sum_{i+j=q}\phi_p(\gamma_i(a)\phi_j(sa)))
=\begin{cases}
\phi_p^2(sa) & \text{ if } p=q\\
0 & \text{ otherwise }
\end{cases}
\]
since $\phi_p$ is additive, and, again, kills elements of the form $ax$.
}\end{remark}

There is obviously a universal $DP$ derivation out of $A$, which we write 
\[
d:A\to\Omega^{DP}_{A/R}\,.
\]
One construction is as follows. Any $R$-module map $s:A\to M$ extends to
a $U(A)$-module map $U(A)\tensor A\to M$. 
As a left $U(A)$-module,
\[
\Omega^{DP}_{A/R}=(U(A)\tensor A)/S
\]
where $S$ is the sub $U(A)$-module generated by the elements
\begin{gather*}
a\tensor b-1\tensor ab+b\tensor a\\
1\tensor\gamma_na-\phi_n\tensor a+\sum_{i+j=n}\gamma_i(a)\phi_j\tensor a\,.
\end{gather*}
In this model, the universal derivation is given by $da=[1\tensor a]$. 

This expression obscures the simplicity of $\Omega^{DP}_{A/R}$ in general. 

\begin{theorem} \label{thm-main}
Let $A$ be a $DP$ algebra. There is a unique $DP$ 
$A$-module 
structure on $\Omega^{CA}_{A/R}$ such that $d:A\to\Omega^{CA}_{A/R}$ is a 
$DP$ derivation, and it serves as the universal $DP$ derivation.
\end{theorem}

The case in which $R$ 
is a field of characteristic $p$ is proven in \cite{dokas-differentials}.

\begin{proof}

We follow \cite{dokas-differentials} in a proof of this theorem. 
To begin with,
we employ one of the standard constructions of the K\"ahler differentials
to construct a $DP$ $A$-module structure on $\Omega^{CA}_{A/R}$. 

The ``fold'' map $\nabla:A\textstyle{\coprod}A\to A$ is characterized 
as a $DP$ algebra map by
the equations $\nabla\circ\inj_1=1_A=\nabla\circ\inj_2$.
The only algebra map satisfying them sends each of the the factors in 
$A\textstyle{\coprod}A=(A\tensor R)\oplus(A\tensor A)\oplus(R\tensor A)$
to $A$ by the product, so the product map is a $DP$ algebra map.

The kernel $I$ of this map is thus a $DP$ subalgebra
of $A\textstyle{\coprod}A$. One construction of $\Omega^{CA}_{A/R}$ is as
the quotient $I/I^2$, which thus has a natural structure of an abelian
$DP$ algebra.

$\Omega^{CA}_{A/R}$ is of course also a module over $A_+$. We have to see 
that this $A_+$-module structure coheres
with the divided powers; that is, that $\gamma_p((a\tensor1)\omega)\in I^2$
for any $\omega\in I$. By the product formula this is 
$(a^p\tensor1)\gamma_p\omega$. But $a^p=p!\gamma_p(a)$ and 
$p\gamma_p\omega \in I^2$ since $p\phi_p=0$ in $I/I^2$.

So $\Omega^{CA}_{A/R}$ is a left $U(A)$-module. 

Next, the universal derivation $d:A\to\Omega^{CA}_{A/R}$ is in 
fact a $DP$ derivation. We refer to \cite{dokas-differentials} for this;
we observe that the argument for Proposition 2.8 in that paper 
does not require one to work over a field of characteristic $p$, but
rather is a mod $p$ result and so is applicable since $p\gamma_p=0$.

Finally, suppose that $s:A\to M$ is any $DP$ derivation. Since it is in
particular a derivation in $\bCA_R$, there is a unique $A_+$-module map
$f:\Omega^{CA}_{A/R}\to M$ such that $s=f\circ d$. We claim that $f$ is in fact
a map of $DP$ $A$-modules.

The $DP$ algebra over $A$ determined by $M$, $A\oplus M$, admits two $DP$
algebra maps from $A$: the unit map and the map corresponding to the 
derivation $s$. Together they induce a $DP$ algebra map 
$A\coprod A\to A\oplus M$. We claim that its restriction to $I$ factors 
through the inclusion $M\hookrightarrow A\oplus M$. To see this let
$\sum a_i\tensor b_i\in I$, so that $\sum a_ib_i=0$ in $A$. Its image in 
$A\oplus M$ is 
\[
\sum(a_i,0)(b_i,s(b_i))=\sum(a_ib_i,a_is(b_i))=(0,\sum a_is(b_i))\,.
\]
Since $M$ is an abelian algebra, products vanish in it, so the map
$I\to M$ factors through the quotient $I/I^2=\Omega^{CA}_{A/R}$. 
This construction of $f:\Omega^{CA}_{A/R}\to M$ makes it clear that it
is indeed a $DP$ module homomorphism. 
\end{proof}

\section{Examples}

There is one case in which the module of K\"ahler differentials is easy
to compute for formal reasons. 

\begin{lemma}[{cf. {\cite[Theorem 2.9]{dokas-differentials}}}]
\label{lemma-free} 
Let $V$ be an $R$-module and let $A=\Gamma_R(V)$ be the free $DP$ 
algebra generated by the $R$-module $V$. 
Then there is a unique $DP$ derivation $d:A\to U(A)\tensor V$
such that $dv=1\tensor v$ for $v\in V$,
and it is the universal $DP$ derivation; 
\[
\Omega^{DP}_{A/R}=U(A)\tensor V\,.
\]
\end{lemma}
\begin{proof}
We simply verify that $(U(A)\tensor V,d)$ satisfies the universal property.
Let $M$ be any $U(A)$-module and consider the corresponding abelian
object over $A$. Let $i:V\to A$ denote the inclusion of generators. Then:
\begin{gather*}
\Der^{DP}_R(A,M)=
\left\{\rule[-3pt]{0pt}{35pt}\right.\raisebox{20pt}
{\xymatrix{& A\oplus M \ar[d] \\
A \ar@{.>}[ur] \ar[r]^1 & A }}
\left.\rule[-3pt]{0pt}{35pt}\right\}_{\bDPA_R}=
\left\{\rule[-3pt]{0pt}{35pt}\right.\raisebox{20pt}
{\xymatrix{& A\oplus M \ar[d] \\
V \ar@{.>}[ur] \ar[r]^i & A }}
\left.\rule[-3pt]{0pt}{35pt}\right\}_{\bMod_R}\\
=\Hom_R(V,M)=\Hom_{U(A)}(U(A)\tensor V,M)\,.
\end{gather*}
We leave the check that the universal differential is as stated to the reader.
\end{proof}

Theorem \ref{thm-main} has the following consequence.

\begin{corollary}[\cite{roby-differentials}]
The module of K\"ahler differentials of the free $DP$ algebra $A=\Gamma_R(V)$ 
when regarded as merely a commutative algebra is
\[
\Omega^{CA}_{A/R}=U(A)\tensor V
\]
as an $A_+$-module. 
\end{corollary}

\begin{example} {\em For example suppose that $V$ is free of rank 1 over $R$,
with generator $x$. Then $A=\Gamma_R(R)$ is the
free $R$-module on $\{\gamma_n(x):n>0\}$. Its universal enveloping algebra has
the form 
\[
A_+\oplus\bigoplus_p(A_+/p)[\phi_p(dx)]
\]
as an $A_+$-module. To lighten notation, abbreviate $\gamma_n(x)$ to 
$\gamma_n$ and $\phi_p(dx)$ as $\phi_p$. 

Now $\Omega^{CA}_{A/R}$ is generated as an
$A_+$-module by the elements $d\gamma_n(x)$ (which we will abbreviate to
$d\gamma_n$). This expression does not reveal its 
structure as an $A_+$-module. Since $d:A\to\Omega_{A/R}$ is a $DP$ derivation, 
we have the relations
\[
d\gamma_n=\phi_n+\gamma_1\phi_{n-1}+\cdots+\gamma_{n-1}\phi_1\,.
\]
The $d\gamma_n$s are determined by the $\phi_j$'s by means of an invertible
lower triangular band matrix. These equations have a unique solution, namely
\[
\phi_n=d\gamma_n-\gamma_1d\gamma_{n-1}+\cdots+(-1)^{n-1}\gamma_{n-1}d\gamma_1\,.
\]
To see this, substitute these 
values for the $\phi_j$'s into the right hand side of the equation.
It is the sum of the following sums:
\begin{alignat*}{5}
d\gamma_n{-}&\gamma_1d\gamma_{n-1}{+}&\gamma_2d\gamma_{n-2}{-}&\cdots
{+}&(-1)^{n-2}\gamma_{n-2}d\gamma_2{+}&(-1)^{n-1}\gamma_{n-1}d\gamma_1\\
&\gamma_1d\gamma_{n-1}{-}&\gamma_1\gamma_1d\gamma_{n-2}{+}&\cdots
{+}&(-1)^{n-3}\gamma_1\gamma_{n-3}d\gamma_2{+}&
(-1)^{n-2}\gamma_1\gamma_{n-2}d\gamma_1\\
&&\gamma_2d\gamma_{n-2}-&\cdots
+&(-1)^{n-4}\gamma_2\gamma_{n-4}d\gamma_2+&(-1)^{n-3}\gamma_2\gamma_{n-3}d\gamma_1\\
&&&&\cdots&\cdots\\
&&&&\gamma_{n-2}d\gamma_2-&\gamma_{n-2}\gamma_1d\gamma_1\\
&&&&&\gamma_{n-1}d\gamma_1
\end{alignat*}
For $k>0$, the $k$th column is of the form
\begin{gather*}
(-1)^k\left(\gamma_k+\sum_{i=1}^{k-1}(-1)^i\gamma_i\gamma_{k-i}+(-1)^k\gamma_k\right)
d\gamma_{n-k}\\
=(-1)^k\left(\sum_{i=0}^k(-1)^i\binom{k}{i}\right)\gamma_kd\gamma_{n-k}\,,
\end{gather*}
since $\gamma_i\gamma_{k-i}=\binom{k}{i}\gamma_k$.
But the alternating sum of binomial coefficients vanishes, leaving only
$d\gamma_n$ as claimed. 

This provides an explict description of $\Omega_{A/R}$ as an 
$R$-module: 
\[
\Omega_{A/R}=A_+\langle d\gamma_1\rangle\oplus\bigoplus_{p}\bigoplus_{e\geq1}
(A_+/p)\langle\phi_{p^e}\rangle\,.
\]
}\end{example}

\begin{remark}{\em  This example actually provides a general expression for
$\phi_n(a)$ in terms of $\{d\gamma_i(a):i\leq n\}$, 
for any $DP$ algebra $A$ and any 
$a\in A$, since there is a unique $DP$ algebra map 
$\Gamma\langle x\rangle\to A$ sending $x$ to $a$. By naturality of the
operators $\gamma_n$ and $\phi_n$, 
\[
\phi_n(da)=d\gamma_n(a)+\sum_{i+j=n}(-1)^i\gamma_i(a)d\gamma_j(a)
\]
in the $U(A)$-module $\Omega_{A/R}$. These elements generate $\Omega_{A/R}$ as
an $A_+$-module. Their dependence on $a$ is additive and ``Frobenius linear'' 
in the sense that $\phi_n(r\,da)=r^n\phi_n(da)$ for $r\in R$. 
And they satisfy the other properties
of the $\phi_n$'s: $\phi_n(da)=0$ unless $n$ is a prime power and  
$p\phi_{p^e}(da)=0$ for any prime number $p$ and any $e\geq1$.
Moreover, the equations in \ref{rem-implications} show that 
$\phi_n(a)$ depends only on the class of $a$ in the ``module of $DP$ 
indecomposables'' of $A$: the maximal quotient $QA$ of $A$ in which products
and divided powers $\gamma_n$ for $n>1$ vanish. So these elements are 
determined by their values on a choice of lifts of $R$-module generators
of $QA$.
}\end{remark}

As observed by Roby, Theorem \ref{thm-main} 
and Lemma \ref{lemma-free} together 
determine the indecomposables in a free divided power algebra:

\begin{corollary}[\cite{roby-differentials}]
The indecomposable quotient of the free $DP$ algebra $A=\Gamma_R(V)$ is
\[
A/A^2=U(0)\tensor V\,.
\]
\end{corollary}
\begin{proof}
First, for any algebra $A$, the augmentation $\epsilon:A_+\to R$ puts an
$A_+$-module structure such that $AV=0$ on any $R$-module $V$. A derivation
$s:A\to V$ satisfies $s(ab)=as(b)+bs(a)=0$, and so factors uniquely through an
$R$-module map $A/A^2\to V$. Conversely, for $R$-module map $A/A^2\to V$
the composite $A\to A/A^2\to V$ is a derivation. This implies
\[
R\tensor_{A_+}\Omega_{A/R}=A/A^2\,.
\]

For $A=\Gamma_R(V)$, we can now calculate
\[
A/A^2=R\tensor_{A_+}\Omega_{A/R}=R\tensor_{A_+}(U(A)\tensor V)=U(0)\tensor V
\]
where the last equality follows from $U(A)=A_+\tensor U(0)$.
\end{proof}

\begin{remark}{\em 
The $DP$ module $\Omega^{DP}_{A/R}$ is the abelianization
of the identity map $A\downarrow A$ as an object of $\bDPA/A$. This special
case of the abelianization functor
\[
\Ab_A:\bDPA/A\to\Ab(\bDPA/A)=\bMod_{U(A)}
\]
-- the left adjoint of the inclusion -- in fact determines the whole functor:
Given $B\downarrow A$,
\[
\Ab_A(B)=U(A)\tensor_{U(B)}\Omega^{DP}_{B/R}\,.
\]
This is easily seen \cite{dokas-triple} using the fact that the pullback
of $A\oplus M\downarrow A$ along $B\to A$ is $B\oplus M|_B\downarrow B$. 
}\end{remark}

\begin{remark} {\em The canonical expression for $\Omega^{CA}_{A/R}$ is
\[
\Omega^{CA}_{A/R}=A_+\tensor A/(a\tensor b-1\tensor ab+b\tensor a)
\]
and the natural map $\Omega^{CA}_{A/R}\to\Omega^{DP}_{A/R}$ includes the 
unit in $U(0)$ and collapses the relation involving the divided powers.
So the theorem shows that the effect of those relations is precisely to
kill the augmentation ideal in $U(0)$. The relations clearly have that
effect; what is in question is whether they do any more. 
}\end{remark}

\end{document}